\documentclass[10pt]{article}

\usepackage{amsmath,amssymb,amsthm}
\def\tr{{\rm tr\,}} \def\Z{\mathbb{Z}} \def\N{\mathbb{N}}

\def\R{\mathbb{R}}\def\C{\mathbb{C}}\def\D{\mathbb{D}}\def\E{\mathbb{E}}\def\Q{\mathbb{Q}}
\title{Random walks in the hyperbolic plane and the Minkowski question mark function}
\author{G\'{e}rard Letac\thanks{Laboratoire de Statistique
et Probabilit\'es, Universit\'e  de Toulouse, Toulouse, France. \texttt{gerard.letac@math.univ-toulouse.fr}},   Mauro
Piccioni\thanks{Dipartimento di Matematica, Sapienza Universit\`{a} di Roma, 00185 Rome, Italia.}\ .}

=-0.5 in
\begin{document}
\maketitle 

\begin{abstract} Consider $G=SL_2(\Z)/\{\pm I\}$ acting on the complex upper half plane $H$ by $h_M(z)=\frac{az+b}{cz+d},$ for $M \in G$. Let $D=\{z \in H: |z|\geq 1, |\Re(z)|\leq 1/2\}$. We consider the set $\mathcal{E} \subset G$ with the $9$ elements $M$, different from the identity, such that $\tr (MM^T)\leq 3$. We equip the tiling of $H$ defined by $\D=\{h_M(D), M \in G\}$ with a graph structure where the neighbours are defined by $h_M(D) \cap h_{M'}(D) \neq \emptyset$, equivalently $M^{-1}M' \in \mathcal{E}$.

The present paper studies several Markov chains related to the above structure. We show that the simple random walk on the above graph converges a.s. to a point $X$ of the real line with the same distribution of $S_2 W^{S_1}$, where
$S_1,S_2,W$ are independent with $\Pr (S_i=\pm 1)=1/2$ and where $W$ is valued in $(0,1)$ with distribution  $\Pr(W<w)=\textbf{?}(w)$. Here $\textbf{?}$ is the Minkowski function. If  $K_1, K_2, \ldots$ are i.i.d with distribution $\Pr (K_i=n)= 1/2^n$ for $n=1,2,\ldots$, then $W= \frac{1}{K_1+\frac {1}{K_2+\ldots}}$: this known result (Isola (2014)) is derived again here.

\vspace{4mm}\noindent

\vspace{2mm}

\vspace{4mm}\noindent \textsc{Acknowledgements} We are grateful to  an anonymous referee for useful  comments and references. G.L. thanks Sapienza Universit\`a di Roma  for its generous support during the preparation of this paper.

\end{abstract}
\newpage
\section{Introduction} In this paper we are concerned with random walks in the upper half of the complex plane (hyperbolic plane)
$$
H=\{z \in \C: \Im (z)>0\}.
$$
The random walk is induced by the action on $H$ of the modular group $G=PSL(2,\mathbb {Z})$, which is the multiplicative group of all $ 2 \times 2$ matrices with integer entries and determinant $1$, quotiented by identifying two matrices when one is equal to the other multiplied by $-I_2$. The basic properties of this action are presented in Serre (1977). For 
\begin{equation}M=\pm \left[\begin{array}{cc}a&b\\c&d\end{array}\right]\end{equation}
with $a,b,c,d\in \Z$ and $\det M=1$ we let
\begin{equation}\label{homographies}
h_M(z)=\frac{az+b}{cz+d},\ \  z \in H.
\end{equation}
It is readily verified that $h_M$ maps $H$ into $H$ since $\Im (h_M(z))$ has the same sign as $\Im (z)$, for any $M \in G$. Morevover $$h_{MM'}(z)=h_M(h_{M'}(z)),$$ which proves that $G$ acts on $H$. Furthermore this action is faithful, as a consequence of the quotient made above.

The discrete nature of $G$ implies the existence of fundamental domains for the action of $G$ on $H$. Roughly speaking, a fundamental domain is a subset of $H$ which contains one element for each $G$-orbit. The traditional choice of a fundamental domain for the above action is the region $D$ of all $z$'s in $H$ such that $|z|\geq 1$ and $|\Re(z)|\leq 1/2$. Each point in the interior of $D$ (and the point $i$) belongs to a different orbit, whereas the remaining orbits intersect two different points on the boundary of $D$. As a consequence 
the set of all images $\D=\{h_M(D), M \in G\}$ covers the whole space, with two distinct sets overlapping at most on their boundaries. As long as $z \in D$ does not lie on the boundary of a set in $\D$, there is a natural "projection" of $H$ onto $\D$. The image $h_M(D)$ intersects $D$ only for $M=\pm E_i$, for some  $i=0,1,\ldots,8$, where

$$\begin{array}{cccc}
E_1=\left[\begin{array}{cc}1&1\\0&1\end{array}\right]&E_2=\left[\begin{array}{cc}1&-1\\0&1\end{array}\right]&E_3=\left[\begin{array}{cc}1&-1\\1&0\end{array}\right]&E_4=\left[\begin{array}{cc}-1&-1\\1&0\end{array}\right]\\&&&\\

h_1(z)=1+z&h_2(z)=-1+z&h_3(z)=1-\frac{1}{z}&h_4(z)=-1-\frac{1}{z}
\end{array}$$

$$\begin{array}{cccc}
E_5=\left[\begin{array}{cc}1&0\\1&1\end{array}\right]&E_6=\left[\begin{array}{cc}1&0\\-1&1\end{array}\right]&E_7=\left[\begin{array}{cc}0&-1\\1&1\end{array}\right]&E_8=\left[\begin{array}{cc}0&-1\\1&-1\end{array}\right]\\&&&\\

h_5(z)=\frac{z}{1+z}&h_6(z)=\frac{z}{1-z}&h_7(z)=-\frac{1}{1+z}&h_8(z)=\frac{1}{1-z}
\end{array}$$

$$
E_0=\left[\begin{array}{cc}0&-1\\1&0\end{array}\right], \\ h_0(z)=-\frac {1}{z},
$$
with $h_i=h_{E_i}$, for $i=0,1,\ldots,8$: see the well known picture at page 128 of Serre (1977). These matrices (and their sign change) are characterized in $G$ as having the trace of $MM^T$ less or equal to three, aside from the identity.

These matrices generate the whole group $G$, since already $\pm E_0$ and $\pm E_1$ together have this property. More generally $h_M(D)$ intersects $h_{M'}(D)$ if and only if $M'=\pm ME_i$ for some   $i=0,1,\ldots,8$. As a consequence one can endow $\D$ with a graph structure, joining intersecting regions with an edge, getting in this way a regular graph with vertices in $\D$ of degree $9$ (incidentally, this graph is isomorphic to the Cayley graph of $G$, taking the $\pm E_i$'s as a set of generators). 

In the following we are going to investigate the asymptotic behaviour of the $H$-valued processes 
\begin{equation}\label{complex}
Z_n^z=h_{M_n} \circ \cdots \circ h_{M_1}(z),V_n^z=h_{M_1} \circ \cdots \circ h_{M_n}(z), z \in H,
\end{equation}
where $(M_k, k=1,\ldots)$ is a sequence of i.i.d. random matrices, taking each of the $9$ possible values $E_i, i=0,1,\ldots,8$, with the same probability $1/9$ (from now on we take for granted the identification of $M$ with $-M$). For each $z \in H$ and each $n=1,2,\ldots$, the complex random variables  $Z_n^z$ and $V_n^z$ have the same law, but whereas $(Z_n^z)_{n=0}^{\infty}$ is  always a homogeneous Markov chain, the process $(V_n^z)_{n=0}^{\infty}$ is not. More specifically, in most of the cases it remains a Markov chain, but with transition probabilities depending on $z$. In fact, if $z$ and $h_{M}(z)$ are known, one deduces $M$, except when the stabilizer of $z$ is non trivial, which happens only for a denumerable subset of $z$ in $H$ (Serre, page 129). Whenever $z \in  h_M(\mathrm{int} D)$, for some $M \in G$, both processes can be "projected" on the vertices of the graph. Of the greatest interest is the fact that, for $z$ in the interior of $D$, $V_n^z$ "projected" on the graph induces a simple nearest neighbour random walk on it.

In order to discuss the asymptotic behaviour of these processes it is necessary to extend the action of $G$ on the boundary $\partial H=\R\cup \{\infty\}$, the extended reals, and to study the processes in $\partial H$
\begin{equation}\label{real}
X_n^x=h_{M_n} \circ \cdots \circ h_{M_1}(x),Y_n^x=h_{M_1} \circ \cdots \circ h_{M_n}(x), x \in \partial H.
\end{equation}
The following result can be found in the second chapter of the book by Bougerol, Lacroix (1985) (see Benoist, Quint (2016) for a more recent presentation). It refers to random compositions of Moebius transformations of the form (\ref{homographies}) as they appear in (\ref{complex}) and (\ref{real}), with general unimodular matrices with real coefficients.

\vspace{4mm}\noindent\textbf{Theorem 1.1} Let $(M_k, k=1,\ldots)$ be an i.i.d. sequence of $2 \times 2$ unimodular matrices and let $\mathcal {G}$ be the smallest closed subgroup which contains the support of their law. Suppose that:
\begin{enumerate}
\item $\mathcal {G}$ is not compact;
\item There does not exist a subset $L$ in $\R^2$ which is a finite union of one-dimensional subspaces which is invariant under all matrices in $\mathcal {G}$.
\end{enumerate}
Then the following hold:
\begin{enumerate}
\item For any $z \in H$, with probability $1$, $(V_n^z)$, as defined in (\ref{complex}), converges to a random variable $Z \in \partial H$, as $n \to \infty$;
\item With probability $1$, for any bounded and continuous function $f$ defined on $\partial H$, $\int f(Y_n^x)\nu (dx)$ converges a.s. to $f(Z)$ as $n \to \infty$, for any probability measure $\nu$ on $\partial H$;
\item The law $\lambda$ of $Z$ is the unique stationary measure for the Markov chain $(X_n^x)$ on $\partial H$, and it is atomless.
\end{enumerate}

The assumptions of the theorem are clearly satisfied when the law of $M$ is supported by the nine values $E_i$, $i=0,1,\ldots,8$. The different kind of convergence stated by  Theorem 1.1 in the complex and in the real case is due to the fact that the product of matrices $M_1 \circ \cdots \circ M_n$, properly normalized, converges to a matrix of rank one, that has a non trivial null space which has is avoided w.p. $1$ by the "initial" vector $(x,1)^t$ since the distribution of this one-dimensional null space is atomless (Corollary 4.8  in Benoist, Quint (2016)).

Thus, from the identity in law of $(Z_n^z)$ with $(V_n^z)$ and of $(X_n^x)$ and $(Y_n^x)$, and the fact that convergence a.s. implies convergence in law, one can deduce the following corollary:

\vspace{4mm}\noindent\textbf{Corollary 1.2} Under the assumption of the previous theorem:
\begin{enumerate}
\item For any $z \in H$, $(Z_n^z)$, as defined in (\ref{complex}), converges weakly to $\lambda$ , as $n \to \infty$;
\item For any atomless probability measure $\nu$ on $\partial H$, the process $(X_n^x)$ as defined in (\ref{real}), with $x$ taken to be $\nu$-distributed, converges weakly to $\lambda$ as $n \to \infty$.
\end{enumerate}

In the present paper the main goal is to identify the unique stationary distribution $\lambda$ for the chain $(X_n^x)$. In order to achieve this goal, we start in Section 2 with the observation that the transition kernel of the chain is equivariant under the action of the four elements group $\Gamma$ on $\R \cup \{\infty\}$, generated by the mappings $g_0(x)=h_0(x)=-1/x$ and $g_1(x)=-x$. This has the consequence that initial laws which are invariant under this group keep this property with the iterations of the Markov chain. Moreover for any function $C$ which is constant on the orbits of this group one obtains that $(C(X_n^x))$  is by itself a Markov chain. By choosing $C(x)=\min \{|x|,\frac {1}{|x|}\}$ we project the dynamics of $(X_n^x)$ from the extended reals to the unit interval $[0,1]$ and we characterize the stationary distribution of this projected Markov chain with the two properties of symmetry w.r.t. $1/2$ and invariance under a certain "tent" map of the interval. In Section 3 we reformulate these invariances in terms of continued fraction expansions, leading to identify the stationary distribution function for $(C(X^x_n))$ as the Minkowski's question mark function $\textbf{?}$ (Minkowski (1904)). A definition of $\textbf{?}$ can be found in (\ref{qm}) below.  The paper by Chassaing, Letac, Mora (1984) can be also consulted for the links between $\textbf{?}$ and the sequences of Farey-Brocot. This function  is a remarkable example of a continuous  singular distribution function on $[0,1]$, of which we are going to review some of its properties. By "lifting" this law on the extended reals to enforce the desired invariance under $\Gamma$, the unique stationary distribution $\lambda$ for $(X_n^x)$ is finally obtained. It turns out that its survival function is a symmetrized version of the so-called Denjoy-Minkowski function of parameter $1/2$ (see Denjoy (1938)). Finally we have to mention that the interest for the Minkowski function $\textbf{?}$ has been recently revived by the proof by Jordan and Sahlsten (2015) of the 1943 Salem conjecture $\lim_n\int_0^1e^{i2\pi nx}d\textbf{?}(x)=0.$

\section{Group invariance properties of the Markov chain $(X_n)$}

The first observation we are going to perform concerns a certain equivariance property of the Markov chain $(X_n^x)$ defined in (\ref{real}) and its consequences. We will use extensively the notation $X \sim \alpha$ when $X$ has the distribution $\alpha$ and $X \sim Y$  when two random variables $X$ and $Y$ have the same law.

\vspace{4mm}\noindent\textbf{Lemma 2.1} Let $g_0(x)=-\frac {1}{x}$ and $g_1(x)=-x$ for $x \in  \R \cup \{\infty\}$.  Let $M_1$ be uniformly distributed on the set of matrices $E_i$, for $i=0,1,\ldots,8$. Then 
\begin{equation}\label{numeroter}
h_{M_1}(g_j(x))\sim g_j(h_{M_1}(x)), j=0,1,
\end{equation}

\begin{proof} Notice that
$$
 h_0(x)=-h_0(-x),\ h_{2i-1}(x)=-h_{2i}(-x), i=1,\ldots,4,
$$
which shows (\ref{numeroter})  with $j=1.$ Similarly
$$
 h_0(1/x)=-\frac {1}{h_0(x)}, \ h_i(-1/x)=-\frac{1}{h_{\varphi(i)}(x)}, i=1,2,5,6.
$$
with $\varphi(1)=3,\ \varphi(2)=4,\ \varphi(5)=8,\ \varphi(6)=7$. From this, formula (\ref{numeroter}) with $j=0$ is obtained.
\end{proof}

It is clear that the statement of the previous lemma holds for any $g$ belonging to the group $\Gamma$ generated by $g_0$ and $g_1$ (shortly, also for $g(x)=1/x$). The consequences of the previous lemma are important.

\vspace{4mm}\noindent\textbf{Proposition 2.2} Let $X_0$ be a random variable on the extended reals with the property $X_0\sim g(X_0)$, for any $g \in \Gamma$. Define $X_1=h_{M_1}(X_0)$, where $M_1$ assumes the values $E_i$ with probability $1/9$, for $i=0,1, \ldots,8$, independently of $X_0$. Then $X_1 \sim g(X_1)$. Furthermore, if $f$ is a bounded measurable function on $\R \cup \{\infty\}$ such that $f(x)=f(g(x))$ for any $g \in \Gamma$, then the function
$s(x)=\E(f(h_{M_1}(x)))$ has again the property $s(x)=s(g(x))$, for any $x \in \R \cup \{\infty\}$.

\begin{proof} For the first statement it is enough to notice that  for any $g \in \Gamma$
$$
g(X_1)=g(h_{M_1}(X_0))\sim h_{M_1}(g(X_0))\sim h_{M_1}(X_0)=X_1.
$$
As far as the second is concerned notice that similarly
$$
f(h_{M_1}(g(x)))\sim f(g(h_{M_1}(x)))=f(h_{M_1}(x)).
$$
\end{proof}

The last statement in the above proposition suggests to introduce a function on the extended reals $\R \cup \{\infty\}$, whose values distinguish among the orbits of $\Gamma$. A convenient function with this property is the function $C:  \R \cup \{\infty\} \rightarrow [0,1]$ defined by
$$
C(x)=\min\{|x|,\frac{1}{|x|}\}.
$$
Given a random variable $X$ with extended real values the distribution of $C(X)$ can be immediately computed. When $X \sim -X$ and $X\sim -1/X$, this transformation can be easily inverted, "lifting" the law of $C(X)$ to  $\R \cup \{\infty\}$, as stated in the next lemma.

\vspace{4mm}\noindent\textbf{Lemma 2.3:} Let $X$ be an extended real valued random variable. Then $X\sim g(X)$ for any $g \in \Gamma$ if and only if the conditional distribution of $X$ given $W=C(X)$ is $\mu_W$, where 
\begin{equation}\label{CMEASURE}\mu_w=\frac{1}{4}\left(\delta_w+\delta_{-w}+\delta_{1/w}+\delta_{-1/w}\right), \ w \in [0,1].\end{equation} 

In this case 
\begin{equation}\label{TRUC}
X \sim S_2 W^{S_1},
\end{equation}
where $S_i, i=1,2$ are independent random variables with $\Pr(S_i=\pm 1)=1/2$, independent of $W$.


\begin{proof} First observe that if $f$ is any bounded measurable function on the extended reals, $x \in \R \cup \{\infty\}$ and $w=C(x)$ one has
$$
\frac {1}{4}(f(x)+f(-x)+f(-1/x)+f(1/x))=\frac {1}{4}(f(w)+f(-w)+f(-1/w)+f(1/w)).
$$
Next observe that $X\sim g(X)$ for any $g \in \Gamma$ if and only if for any $f$ as above and any $a$ bounded measurable function on the unit interval $[0,1]$, it holds
\begin{eqnarray*}
\E(f(X)a(W))&=&\frac {1}{4}(\E(f(X)a(W)+f(-X)a(W)+f(-1/X)a(W)+f(1/X)a(W)))\\&=&\frac {1}{4}(\E(f(W)a(W)+f(-W)a(W))+f(-1/W)a(W))+f(1/W)a(W)))
\end{eqnarray*}
which yields the first statement of the lemma. For the last distributional representation, it is clear that the r.h.s. of (\ref{TRUC}) satisfies the assignment of the conditional distributions (\ref{CMEASURE}).
\end{proof}


To continue, define the following mappings of the unit interval $[0,1]$ into itself, namely 
\begin{equation}\label{NPT0}
H_0(x)=x, H_1(x)=\frac{1}{1+x},\ H_2(x)=1-x,\ H_3(x)=\min\{\frac{x}{1-x},\ \frac{1-x}{x}\}, H_4(x)=\frac{x}{1+x}.
\end{equation}
It is immediately verified that 
\begin{equation}\label{CPROCESS1}
C(h_i(x))=C(H_i(x)), \,\,\,i=0,1,\ldots,4,
\end{equation}
\begin{equation}\label{CPROCESS2}
C(H_i(x))=C(h_{\psi(i)}(x)),\ \psi(1)=7,\ \psi (2)=8,\ \psi (3)=6,\ \psi (4)=5.
\end{equation}
Thus, given the sequence $(M_n)$ of independent random matrices with uniform distribution on $E_i, i=0,1,\ldots,8$, define $I_n=0$ for $M_n=E_0$ and $I_n=i$ for $M_n=E_i$ or $M_n=E_{\psi(i)}$, $i=1,2,3,4$. One obtains an i.i.d. sequence $(I_n)$ with values in $\{0,1,2,3,4\}$ with distribution $\rho$ such that $\rho(\{0\})=1/9$ and $\rho(\{i\})=2/9$, for $i=1,2,3,4$.

We are now ready to prove the following

\vspace{4mm}\noindent\textbf{Proposition 2.4:} Let $(X_n^x)$ be defined in (\ref{real}). The process $(W_n^w=C(X_n^x))$, with $w=C(x)$, is a Markov chain with values in the unit interval $[0,1]$ which evolves in the following way

\begin{equation}W_{n+1}^w=H_{I_{n+1}}(W_{n}^w), n=0,1,\ldots, W_{0}^w=w\in [0,1] \label{NPT}\end{equation} 

Moreover, if $x=X_0$ with $X_0 \sim g_i(X_0)$, for $i=0,1$, then, for any positive integer $n$, the conditional distibution of $X_n^{X_0}$ given $W_n^{W_0}=w$, where $W_0=C(X_0)$, is given by $\mu_w$, as defined in (\ref{CMEASURE}).

\begin{proof} Collecting together the definitions (\ref{NPT0}), the relations (\ref{CPROCESS1}) and (\ref{CPROCESS2}), Proposition 2.2 and Lemma 2.3 the result is readily obtained.
\end{proof}

As a consequence we have the following

\vspace{4mm}\noindent\textbf{Corollary 2.5:} Let $\nu$ be a stationary distribution for the process $(W_n^w)_{n=1}^{\infty}$ defined in (\ref{NPT}) and let $W \sim \nu$. Define the random variable $X$ as in (\ref{TRUC}). Then the distribution of $X$ is stationary for the chain $(X_n^x)_{n=1}^{\infty}$.

\medskip

Thus for each stationary distribution for the process $(W_n^w, n=1,2,\ldots)$ we can construct a corresponding stationary distribution for the original process $(X_n^x, n=1,2,\ldots)$ by the operation of "lifting" described above. Now we reduce the construction of a stationary distribution for the Markov chain $(W_n^w)$ to the existence of a law which is invariant under two transformations of the unit interval, the symmetry transformation $H_2$ around $1/2$ and the tent-like map $H_3$. The basic point is that the inverse graph of the latter is the union of the graphs of $H_1$ and $H_4$. 

\vspace{4mm}\noindent\textbf{Proposition 2.6:} Let $H_i, i=0,1,\ldots,4$ be defined as in (\ref{NPT0}) and let $W_0$ be a random variable with values in $[0,1]$  with the properties
\begin{equation}\label{INVH}W_0 \sim 1-W_0=H_2(W_0), W_0 \sim \min (\frac {W_0}{1-W_0}, \frac {1-W_0}{W_0})=H_3(W_0).\end{equation}
Then $W_1=H_I(W_0) \sim W_0$, $I$ having the law $\rho$, independent of $W_0$.

\begin{proof} It is trivially $W_0 \sim H_0(W_0)$ and by assumption $W_0 \sim H_2(W_0)$ and $W_0 \sim H_3(W_0)$. Therefore the result holds if one proves that the assumptions on the law of $W_0$ imply that $W_0 \sim H_J(W_0)$, $J$ being a random variable assuming the values $1$ and $4$ with the same probability. It is not too complicate to realize that the law of $W_0$ cannot have an atom at $1/2$. This comes from the fact that the invariance of the law under $H_3$ imply that $1$ has an atom with the same weight.  But this is impossible since $0$ and $1$ are both sent to $0$ by $H_3$, which  contradicts the invariance of the distribution of $W_0$ by $H_3.$

To continue the proof of Proposition 6,  we need the following observation, which for later use is collected as a lemma.

\vspace{4mm}\noindent\textbf{Lemma 2.7:} If $W_0 \sim 1-W_0$, $W_0$ not having an atom at $1/2$, then the law of $H_3(W_0)$ conditional to $W_0< 1/2$ coincide with the law of $H_3(W_0)$ conditional to $\{W_0> 1/2\}$. Thus both coincide with the unconditional law of $H_3(W_0)$.

\begin{proof} Being $H_3(w)=H_3(1-w)$ for any $w \in [0,1]$, the law of $H_3(W_0)$ conditional to $W_0< 1/2$ is equal to the law of $H_3(1-W_0)$ conditional to $W_0<1/2$. Next, replace $W_0$ with $1-W_0$, these being equal in law, to obtain that the law of $H_3(W_0)$ conditional to $1-W_0<1/2$ is still the same. The last conditioning being the same as $W_0> 1/2$, the proof of the first statement of the lemma is finished. The second is obtained from the law of total probabilities. 
\end{proof}

\textit{Proof of Proposition 2.6, continued.} Next observe that $H_4(H_3(w))=w$ for $w < 1/2$ and $H_1(H_3(w))=w$ for $w > 1/2$. As a consequence for $w > 1/2$, from $H_3(W_0) \sim W_0$ one gets that
$$
\Pr(H_J(W_0)>w)=\frac {1}{2}\Pr(H_1(W_0)>w)=\frac {1}{2}\Pr(H_1(H_3(W_0))>w)
$$
$$
=\Pr(W_0 > \frac {1}{2}, H_1(H_3(W_0))>w)=\Pr(W_0>w),
$$
and for $w < 1/2$ 
$$
\Pr(H_J(W_0)<w)=\frac {1}{2}\Pr(H_4(W_0)<w)=\frac {1}{2}\Pr(H_4(H_3(W_0))<w)
$$
$$
=\Pr(W_0 < \frac {1}{2}, H_4(H_3(W_0))<w)=\Pr(W_0<w).
$$
These two together easily imply that $H_J (W_0) \sim W_0$.
\end{proof}

\section{Minkowski's question mark function and Denjoy-Minkowski distribution on the real line}

The goal of this section is to deduce from the invariance properties assumed in (\ref{INVH})  a unique law for $W_0$, whose distribution function is the question mark function introduced by Minkowski. This characterization is well known (see Isola (2014), Lemma 4.1), but here we give 
a probabilistic proof of it. For this purpose the continued fraction representation of irrational numbers in the unit interval $[0,1]$ is required (see Olds (1963)). On this interval we define the function $A_k(w)=\frac {1}{k+w}$ for $k\in \N^{+}$. Likewise define for $k_1,\ldots, k_n,\ldots$ in $\N^{+}$ 
\begin{equation}\label{comp}
A_{k_1,\ldots,k_n}(w)=A_{k_1}\circ \ldots \circ A_{k_n} (w)=\frac{1}{k_1+\frac{1}{k_2+\frac{1}{\ddots+\frac{1}{k_n+w}}}}, \,\,\ n \in \N^+.
\end{equation} 
Then 
\begin{equation}\label{cf}
x=\lim_n A_{k_1,\ldots,k_n}(w)\stackrel{\mathrm{def}}{=}[k_1,\ldots,k_n,\ldots]
\end{equation}
always exists and does not depend on $w \in [0,1]$. Such an $x$ is necessarily an irrational number. Conversely for any irrational number $x\in (0,1)$ there exists a unique sequence $(k_n)_{n \in \N^{+}}$ such that (\ref{cf}) holds.
This is called the continued fraction expansion of $x$: its definition implies the recursion
\begin{equation}\label{contfrac}
x=[k_1,k_2, \ldots]=\frac {1}{k_1+[k_2,k_3, \ldots]}.
\end{equation}

The above construction allows to associate to any probability distribution $p$ on the positive integers an atomless law $\mu (p)$ on the interval $[0,1]$ in the following way. Let $(K_n, n=1,2,\ldots)$ be a sequence of i.i.d. $p$-distributed random variables: then $W=[K_1,K_2,\ldots]$ has the law $\mu (p)$. The function $p \mapsto \mu(p)$ is clearly injective, since $K_1$ is the integer part of $W^{-1}$. The distribution $\mu (p)$ can be characterized as the unique stationary distribution for the Markov chain  $(U_n^u, n=1,2,\ldots)$, where
$$
U_{n+1}^u=A_{K_{n+1}}(U_n^u)=\frac {1}{K_{n+1}+U_n^u}, n=0,1,\ldots, U_0^u=u \in [0,1].
$$
This is an instance of a general principle (see Letac (1986) and Chamayou, Letac (1991), Proposition 1). An equivalent way of stating this property is the following: for $W$ and $K$ independent random variables, with values in $[0,1]$ and $\N^+$ , respectively, it holds 
\begin{equation}\label{charact}
K\sim p, W \sim \frac {1}{K+W} \Longrightarrow W \sim \mu(p).
\end{equation}

Now we are in a position to prove the following result. 

\vspace{4mm}\noindent\textbf{Theorem 3.1} Let $W$ have an atomless law on the interval $[0,1]$: Then the following are equivalent: 
\begin{enumerate}
\item 
$W \sim 1-W$ and $W \sim  \min \{\frac {W}{1-W}, \frac {1-W}{W}\}$.
\item
$W \sim \mu (p)$, with
\begin{equation}\label{GEOM}
p(n)=2^{-n}, n=1,2,\ldots.
\end{equation}
\item
 The distribution function of $W$ at irrational points is given by
\begin{equation}\label{qm}
P(W<[k_1,k_2, \ldots])= 2\sum_{n=1}^{\infty}(-1)^{n+1} 2^{-\sum_{j=1}^n k_j} \stackrel{\mathrm{def}}{=} \textbf{?}([k_1,k_2,\ldots])
\end{equation}
for $k_j=1,2,\ldots$, $j=1,2,\ldots$.
\end{enumerate}

The function $\textbf{?}$ defined in (\ref{qm}) on the irrational numbers is called the Minkowski's question mark function. Being continuous, it can be uniquely extended to the whole unit interval. In fact it is known since the work of Salem (1943) that $$|\textbf{?}(x)-\textbf{?}(x')|\leq C|x-x'|^{\alpha}\  \mathrm{with} \ \alpha=\frac {\log 2}{2 \log \theta}$$ where $\theta=\frac {1+ \sqrt {5}}{2}$ is the golden ratio and $C$ is a constant. The function $\textbf{?}$ is strictly increasing but it is singular w.r.t. the Lebesgue measure, since its derivative is zero a.s. (Salem (1943), Viader,  Bibiloni and Paradis (1998)). A direct proof of the characterization 2) of the distribution function $\textbf{?}$ stated in the previous theorem can be found in Isola (2014), Lemma 4.6.

\begin{proof} 1) implies 2). Since $W$ has an atomless law, we can assume that it takes values in the set of irrationals, and write $W=[K_1,K_2,\ldots]$, with the law of the process $(K_n, n=1,2,\ldots)$ to be determined. Next observe that for any $w=[k_1,k_2,\ldots] \in [0,1]\setminus \Q$ one has for $k_1>1$ (i.e. $w<1/2$) and $k_1=1$ (i.e. $w>1/2$), respectively
\begin{equation}\label{RCF}
\frac {w}{1-w}=[k_1-1,k_2,\ldots], \ \ \frac {1-w}{w}=[k_2,k_3,\ldots].
\end{equation}
We now show by induction the following facts

\begin{itemize}
\item
$(A)_n \,\,\,\, \Pr(K_1=k)=\frac {1}{2^k}, k=1,\ldots, n,$
\item
$(B)_n \,\,\,\,  [K_2,K_3.\ldots] | \{K_1=n\} \sim W,$
\item
$(C)_n \,\,\,\, [K_1-n,K_2,\ldots] | \{K_1>n\} \sim W.$
\end{itemize}

For $n=1$ $(A)_1$ is a consequence of the symmetry of the law of $W$ around $1/2$, whereas $(B)_1$ and $(C)_1$ are obtained from Lemma 2.7: indeed there we established that both the law of $H_3(W)=\frac {W}{1-W}$ conditional to $W<\frac {1}{2}$ and the law of $H_3(W)=\frac {1-W}{W}$ conditional to $W>\frac {1}{2}$, are equal to the unconditional law of $H_3(W)$, which in turn is equal to the law of $W$.

Now assume that $(A)_n, (B)_n, (C)_n$ are true and proceed by induction. Since 
$$
\Pr(K_1=n+1)=\Pr(K_1>n)\Pr(K_1=n+1|K_1>n)
$$ 
and the first factor  by the induction assumption $(A)_n$ is equal to $1/2^n$, we have to prove that $\Pr(K_1=n+1|K_1>n)=\Pr(K_1-n=1|K_1-n>0)=\frac {1}{2}$. This is a consequence of $(C)_n$ and $(A)_1$. Hence $(A)_{n+1}$ is proved.

To prove $(B)_{n+1}$ we condition the l.h.s. of $(C)_n$ by $\{K_1=n+1\}$. We get
$$
[K_1-n,K_2,\ldots]|\{K_1=n+1\} \sim W |\{K_1=1\} =[1, K_2, K_3,\ldots]
$$
so that 
$$
[K_2,K_3,\ldots]|\{K_1=n+1\} \sim [K_2,K_3,\ldots]|\{K_1=1\} \sim W
$$ 
from $(B)_1$. Hence $(B)_{n+1}$ is proved. 

Finally in order to prove $(C)_{n+1}$ we condition the l.h.s.  of $(C)_n$ by $\{K_1>n+1\}$. We get 
$$
[K_1-n,K_2,\ldots]|\{K_1>n+1\} \sim W |\{K_1>1\}=[K_1,K_2,\ldots]|\{K_1>1\}
$$
and this in turn implies that
$$
[K_1-n-1,K_2,\ldots]|\{K_1>n+1\} \sim [K_1-1,K_2,\ldots]|\{K_1>1\}\sim W.
$$
 from $(C)_1$. Hence $(C)_{n+1}$ is proved. Finally notice that $(A)_n$ and $(B)_n$, for any positive integer $n$, are equivalent to the l.h.s. of (\ref{charact}), with $p$ given in (\ref{GEOM}). This establishes 2).

2) implies 3). Consider the representation $W \sim [K_1,K_2,\ldots]$, where $(K_n, n=1,2,\ldots)$ is an i.i.d. sequence of random variables with the same distribution (\ref{GEOM}). The survival function of $K_1$ being equal to $P(K_1 \geq k_1)=2\times 2^{-k_1}$, the events 
$$
E_{n}=\{K_1=k_1,\ldots,K_{n-1}=k_{n-1},K_n\geq k_n\}, n=1,2\ldots
$$
have probabilities
\begin{equation}\label{Geome}
\Pr(E_n)=2\times 2^{-\sum_{j=1}^n k_j}.
\end{equation}
Next another sequence $(F_n)$ is constructed by means of the following  recursion, starting from $F_1=E_1$, 
$$
F_{2n}=F_{2n-1}\setminus E_{2n}, F_{2n+1}=F_{2n}\cup E_{2n+1}, n=1,2,\ldots.
$$ 
The fundamental property is that, for any positive integer $n$
$$
F_{2n} \subset \{W<[k_1,k_2, \ldots]\} \subset F_{2n-1},
$$
since the functions $A_{k_1,\ldots,k_n}$ are decreasing for $n$ odd and increasing for $n$ even, and the range of $A_{k_1,\ldots,k_{n-1},k_{n}+1}$ is an interval adjacent to the right (left) to the range of $A_{k_1,\ldots,k_{n-1},k_{n}}$ if $n$ is even (odd).
By the properties of continued fraction expansions, both the sequences $(F_{2n-1})$ and $(F_{2n})$ converge (from above and from below, respectively) to the event $\{W<[k_1,k_2, \ldots]\}$. Since for any positive integer $n$
$$
\Pr(F_{2n})=\Pr(F_{2n-1})-\Pr(E_{2n}), \Pr(F_{2n+1})=\Pr(F_{2n})+\Pr(E_{2n+1})
$$
it suffices to substitute the expressions (\ref{Geome}) to get the desired (\ref{qm}).

3) implies 1). It consists in a simple verification. Since for $k_1>1$ we have
$$
1-[k_1,k_2,\ldots]=[1,k_1-1,k_2.\ldots],
$$
in order to prove that $\textbf{?}$ corresponds to a probability measure which is symmetric around $1/2$, it is enough to verify
$$
\textbf{?}([k_1,k_2,\ldots])+\textbf{?}([1,k_1-1,k_2,\ldots])=1, k_1>1,
$$
which is straightforward. The second invariance property is deduced from (\ref{RCF}) and   from the fact that $K-1|\{K>1\} \sim K$ when $K$ has the distribution (\ref{GEOM}). 
\end{proof}

Next, by Proposition 2.6 and Theorem 3.1 we have the following

\vspace{4mm}\noindent\textbf{Corollary 3.2} The function $\textbf{?}$ is a stationary distribution function for the Markov chain $(W_n^w)$ defined in (\ref{NPT}). 

By Lemma 2.3, for completing our program we need to compute the law of $S_2W^{S_1}$, where $W,S_1,S_2$ are independent,          $W$ has the distribution function $\textbf{?}$ on $[0,1]$ and  $S_1$ and $S_2$ are two random variables which assume the values $-1$ and $+1$ with the same probability $1/2$. This law is stationary for the process $(X_n^x)$ described in (\ref{real}) because of Proposition 2.2 and Corollary 3.2 and  it is unique by Theorem 1.1. 

As a first step we prove that the distribution of $W^{S_1}$ is the so-called Denjoy-Minkowski function $\chi_{1/2}$ of order $1/2$ (Chassaing \textit{et al.} (1984) page 41). In order to define it, we write positive irrational numbers $y$ in the form 
$$
y=[k_0;k_1,k_2,\ldots]=k_0+\frac {1}{k_1+\frac {1}{k_2+\ldots}},
$$
where $k_0=[y]$ and $y-[y]=[k_1,k_2,\ldots]\in (0,1)$. Now define
\begin{equation}\label{ENJOY}
\chi_{1/2}(y)=\chi_{1/2}([k_0;k_1,k_2,\ldots])=\sum_{n=0}^{\infty} (-1)^n 2^{-\sum_{j=0}^nk_j}.
\end{equation}
 As for the function $\textbf{?}$, it is observed that $\chi_{1/2}$ is a continuous function, thus it extends uniquely to the whole non-negative real line.

\vspace{4mm}\noindent\textbf{Proposition 3.3} Let $W$ have the distribution function $\textbf{?}$. Let $Y=W^{S_1}$. Then the survival function $\Pr(Y>y)$ of $Y$ is the function $\chi_{1/2}(y)$. 

\begin{proof} By comparing (\ref{qm}) with (\ref{ENJOY}) it is immediately verified that for $y$ irrational
$$
\Pr (Y>y)=1-\Pr (S_1=1, X<y)=1-\frac {\textbf{?}(y)}{2}, 0<y<1,
$$
and
$$
\Pr (Y>y)= \Pr(S_1=-1, \frac {1}{X}>y)=\frac {1}{2}\textbf{?}(\frac{1}{y}), y>1.
$$
Now it remains to verify that the r.h.s. of the above expressions coincide with $\chi_{1/2}(y)$, for all irrationals $y$. For the former, this is immediately verified. For the latter, we conclude with the observation 
 that, for $k_0\geq 1$ it is \begin{equation}\label{INVK}\frac {1}{[k_0;k_1,k_2,\ldots]}= [0;k_0, k_1, \ldots].\end{equation} 
\end{proof}

Here  is a noteworthy property of $\chi_{1/2}$.

\vspace{4mm}\noindent\textbf{Proposition 3.4} Let $Y$ be a positive random variable with the survival function $\chi_{1/2}$. Then $[Y]+1$ has the geometric distribution (\ref{GEOM}) and it is independent of $Y-[Y]$, which has the distribution function $\textbf{?}$.  In other words, if $K_0,K_1,\ldots$ are i.i.d. with distribution (\ref{GEOM})   then \begin{equation}\label{LAWY}Y\sim[K_0-1;K_1,K_2,\ldots].\end{equation}

\begin{proof} If $W$ have distribution function $\textbf{?}$, we know that one can  construct $W = [K_1,K_2,\ldots]$, where $(K_n)$ is an i.i.d. sequence of random variables with the distribution (\ref{GEOM}). Moreover, let $S_1$ independent of $W$ such that $\Pr(S_1=\pm 1)=1/2.$  From  Proposition 3.3 we write $ Y=W^{S_1}.$ Thus the law of $Y$ is a mixture, with equal weights, of the law of $[0;K_1,K_2,\ldots]$ and, from (\ref{INVK}),  of that of $[K_1;K_2,\ldots]$. From this one obtains (\ref{LAWY}).
\end{proof}
The last step that ends the determination of the unique stationary distribution $\lambda$ of the Markov chain $(X_n^x)$ defined in (\ref{real}), is a simple symmetrization of the Denjoy-Minkowski function.

\vspace{4mm}\noindent\textbf{Proposition 3.5} Let $X \sim \lambda$, the unique stationary distribution of the chain $(X_n^x)$ defined in (\ref{real}). Then, for any $x\geq 0$
$$
\Pr(X>x)=\Pr(X<-x)=\frac {1}{2} \chi_{1/2}(x)
$$

\begin{proof} It is immediately obtained from the representation $X=S_2Y$, where $Y$ has the survival function $\chi_{1/2}$ and $S_2$ is an independent random variable such that $\Pr(S_2=\pm 1)=1/2.$  
\end{proof}

\section{References}
\vspace{4mm} \noindent\textsc{Bougerol, P. and  Lacroix, J.} (1985) \textit{'Products of random matrices with applications to Schrodinger operators'}, Birkhauser, Boston, MA.

\vspace{4mm} \noindent\textsc{Benoist, Y. and Quint, J.F.} (2016) \textit{'Random walks on reductive groups'}, Springer International Publishing.

\vspace{4mm} \noindent\textsc{Chamayou, J.-F. and  Letac, G.} (1991)
'Explicit stationary distributions for composition of random functions and products of random matrices,' \textit{J. Theoret. Probab.} \textbf{4}: 3-36.

\vspace{4mm} \noindent\textsc{Chassaing, P., Letac, G. and  Mora, M.} (1984)
 'Brocot sequences and random walks in $SL(2,\R)$,
 {\it Springer Lectures Notes, Probability on Groups IX} {\bf
1084}: 37-50.

\vspace{4mm} \noindent \textsc{Denjoy, A.} (1938) 'Sur une fonction r\'eelle de Minkowski',
\textit{J. Math. Pures Appl., S\'er. 17} \textbf{IX}: 105-151.

\vspace{4mm} \noindent \textsc{Isola, S.} (2014)  'Continued fractions and dynamics', \textit{Appl. Math.} \textbf{5}: 1067-1090.

\vspace{4mm}\noindent \textsc{Jordan, T.  and Sahlsten, T.} (2016) 'Fourier transforms of Gibbs measures for the Gauss map',   \textit{Math. Ann.} \textbf{364}: 983-1023.

\vspace{4mm}\noindent \textsc{Minkowski, H.} (1904) 'Zur geometrie der Zahlen', \textit{Verhandlungen des III internationalen \newline Mathematiker-Kongress in Heidelberg}, Berlin.

\vspace{4mm} \noindent \textsc{Olds, C.D.} (1963) \textit{'Continued Fractions'},  The Mathematical Association of America, Washington D.C.

\vspace{4mm} \noindent \textsc{Salem, R.} (1943) 'On some singular monotonic functions which are strictly increasing',  \textit{Trans. Amer. Math. Soc.} \textbf{53}: 427-439.

\vspace{4mm} \noindent \textsc{Serre, J.-P.}(1977) \textit{'Cours d'arithm\'etique'},  deuxi\`eme \'edition revue et corrig\'ee, Presses Univ. de France, Paris.

\vspace{4mm} \noindent\textsc{Viader, P., Paradis, J. and  Bibiloni, L.} (1998) 'A new light  on Minkowski's \textbf{?}(x) function', \textit{J. Number Th.} \textbf{73}: 212-227.

\end{document}